\documentclass[12pt]{amsart}
\usepackage{geometry}                
\usepackage{graphicx}
\usepackage{amssymb}
\usepackage{epstopdf}
\DeclareMathOperator*{\esssup}{\mathtt{esssup}}

\newtheorem{thm}{Theorem}[section]

\newtheorem{prop}[thm]{Proposition}

\newtheorem{lemma}[thm]{Lemma}
\newtheorem{assmp}{Assumption}

\def\R{\mathbb{R}}
\def\dd{\mathrm{\:d}}

\begin{document}
\title{American Options in the Hobson-Rogers Model }
\author{Narn-Rueih Shieh }
\date{}
\maketitle

\begin{abstract}
In this article, we consider a risky asset $X$ for which evolution follows a model proposed by D.G. Hobson and L.C.G. Rogers\cite{HR98}.  We assume that the volatility of $X$ depends on the ratio of the present value and
the exponentially weighted average of the past value. Using the Markovian modelling of the enlarged two-dimensional process, we show that,  for   the  American put option with $X$ as the underlying asset, the
continuation region and the stopped region are separated   a  striking curve . This striking  curve lies between the two striking curves from the basic BSM model, yet is {\em not} monotone.

\end{abstract}

\section{Introduction}
In this article, we propose the optimal stopping problem for a diffusion-type process $X(t)$ in $\R $,
which is governed by a certain stochastic delay differential equation (SDDE).
The main result is to present the optimal striking for the  American put option for which the the underlying risky asset $X(t)$
obeys a delay geometric Brownian motion (DGBM) which was proposed  by Hobson and Rogers\cite{HR98}; we use their offset function of order 1 in the below. The model in \cite{HR98} is also regarded with  the stochastic volatility and
with the volatility smile; yet it has the advantage that  the model preserves the market completeness of which the usual SV model is lack.
\par
We remark that  the SDDE has been under active research for years, useful references are \cite{M98} and \cite{EOS00}.
The applications to European options for which underlying stock is DGBM are well studied; we mention \cite{AHMP07, KSW07, MS13}, among others. Optimal stopping problem related to DGBM appeared in \cite{GR06, GR02}. The feature for the delay equation is that the solution brings the memory from the past, and thus the vast literature on the Markovian solution of an SDE is not readily applicable. 

This article is organized as follows. In Section 2, we formulate our SDDE and propose an enlargement of the dimension to fit the Markvian setting in the 1+1 dimension.
In Section 3, we present our main result, namely to consider the underlying risky asset following a certain DGBM, which was proposed in \cite{HR98},
and to present the optimal striking curve for the associated  American put option. The proofs of our results are given in the  Section 4. The final Section 5 is the conclusion, 
in which we discuss the novelty of the result in the article and some related direction in financial economics. 

{\bf Acknowledgement:} The content of this article is  achieved  mainly while the author  visits 
 Department of Applied Mathematics, The Hong Kong Polytechnic University, during the Spring 2014; I thank Prof.  Zuoquan Xu  for the host.

\section{The asset model}
We consider the following SDDE:
\begin{equation}\label{sdde}
\dd X(t)=b(X(t), Y(t))\dd t+\sigma(X(t), Y(t))\dd B(t),\quad t\geqslant 0,
\end{equation}
in which $B(t)$ is the standard Brownian motion. The process $Y$ defined by
$$
Y(t)=\frac{\int_{-\infty}^{0}e^{\lambda s} X(t+s) \dd s}{\int_{-\infty}^{0}e^{\lambda s}  \dd s}=\lambda\int_{-\infty}^{0}e^{\lambda s} X(t+s) \dd s.
$$
The parameter $\lambda>0$ for the exponential averaging, and the continuous deterministic past-memory
$$
X(s)=\xi (s), \quad s\in (-\infty, 0],
$$
are pre-given.
\par
We notice that the $Y$ has the differential
\begin{equation}\label{sddey}
\dd Y(t)=\lambda ( X(t)- Y(t))\dd t.
\end{equation}
\par
Under suitable Lipschitz and growth conditions, there exits a
unique strong solution for \eqref{sdde}; see \cite{M98} for the
detail. In view of \eqref{sdde} and \eqref{sddey}, we see that the
two-dimensional process $(X(t), Y(t))$ constitutes  a Markovian
process in $\mathbb{R}^2$ with continuous paths. We should remark
that, however, the one-dimensional $X(t)$ is not Markovian. The
following {\it It\^{o} formula for delay} is adapted from
\cite{EOS00}. Let $F(x,y)$ be a given twice differential function
in $(x,y)$.
\begin{align*}
\dd F(X(t), Y(t))=&\left(b(X(t), Y(t))\frac{\partial F}{\partial x}+\lambda (X(t)- Y(t))\frac{\partial F}{\partial y}+\frac{1}{2}\sigma^2(X(t), Y(t))\frac{\partial^2 F}{\partial x^2}\right)\dd t\\
&+\left(\sigma(X(t), Y(t))\frac{\partial F}{\partial x}\right)\dd B(t).
\end{align*}

We should remark that, in generality as it is presented in  \cite{M98}, an SDDE  $X(t)$ is regarded as an {\it infinite-dimensional }\/ Markov process. Here we have a two-dimensional Markovian enlargement  $(X(t), Y(t))$  is due to the choice of the memory $Y(t)$. It is seen that $(X(t), Y(t))$ is \emph{not} Markovian if we choose $Y(t)= X(t-t_0)$ for some time-instant $t_0$; this type of delay indeed appears in the  main literature of SDDE, see the final Section 5 for some discussions.  We also remark  that, the above It\^{o} formula  for delay  is so-named is mainly 
to remind the reader for such formula existing in  the general SDDE context; the above one can indeed  be induced from the two-variate It\^{o} formula  for the two-dimensional Markovian diffusion  $(X,Y)$, as shown in   \S 6.6 of \cite{S04}  

\section{Main Result}
In this main section, we consider the  American put option for which the underlying risky asset $X(t)$ is with the constant risk-less interest rate $r>0$, and is with the volatility $\sigma(s)$ depending solely on the ratio of the present-value $X(t)$
 and the past-value $Y(t)$, where $Y$ is given by \eqref{sddey}. Thus, $X(t)$ obeys the following DGBM
\begin{equation}\label{dgbm}
 \dd X(t)=r X(t)\dd t+\sigma(Z(t)) X(t)\dd B(t),\ t\geq 0,
\end{equation}
in which $Z$ is the ratio process
$$Z(t)=\frac{X(t)}{Y(t)}.$$
 By It\^{o} formula for delay in Section 2, we have the following differential for $Z$,
 \begin{equation}\label{sddez}
 \dd Z(t)=(r+\lambda -\lambda Z(t))Z(t)\dd t+\sigma(Z(t)) Z(t)\dd B(t).
 \end{equation}
 We remark that, from \eqref{dgbm} and \eqref{sddez}, we have a certain stochastic volatility model for the risky asset $X$. However, since we have only one source of randomness, namely $B(t)$, the market is still complete; see the Remark 3.2 in \cite{HR98} for this aspect.
\par
The $(X(t), Z(t))$ is a strong Markovian process in $\mathbb{R}^2$ with continuous paths; we also remark that the log-process $\ln Z$  has a certain mean reverting property, indeed by applying It\^{o} formula to $Z$ and
$\ln z$ we have
$$
Z(t)=Z(0)\exp\big\{\int_{0}^{t} [-\lambda (Z(t)-(1+ \frac{r}{\lambda}-\frac{\sigma^2(Z(s))}{2\lambda}))]\\ds + \int_{0}^{t}\sigma(Z(s))\dd B(s)          \big\}.
$$
By the  Assumption \ref{ubd} in the below,   we may see that $\ln Z(t)$ is mean-reverting to a zone $[1+ \frac{r}{\lambda}-\frac{\sigma^2_1}{2\lambda},
 1+ \frac{r}{\lambda}-\frac{\sigma^2_2}{2\lambda}].$
We remark that,  a diffusion which exhibits such mean-reverting to a zone,  rather than to a constant or to a time-curve, seems to be a new class, if we compare with   the usual mean-reverting Ornstein-Uhlenbeck process; a discussion such ``zone-reverting" diffusions will appear in   \cite{Sh18}.

 The fair price of the  American put option in the time-horizon $[0,T]$, associated with the two-dimensional Markovian process $(X(t), Z(t))$, is defined to be
 $$
 V(x,z)=\esssup_{\tau}\mathbb{E}_{(x,z)}\big[e^{-r\tau} (K-X(\tau))^{+}\big],
 $$
where the notation $\mathbb{E}_{(x,z)}[\cdot]$ denotes the expectation w.r.t. the process starting at $(x,z)$. The $\tau$ is ranging over the class of all stopping times over the time-horizon $[0, T]$, w.r.t. the Brownian filtration $\{\mathcal{F}(t), t\geq 0\}$; see, for example, Section 25.1 of \cite{PSh06}.
\par
Our standing assumption is that, besides the continuity of $\sigma (z)$,
\begin{assmp}\label{ubd}
 $0< \sigma_2=\inf \sigma (z)<\sup \sigma (z)=\sigma_1< \infty. $
\end{assmp}
This is a reasonable assumption for the non-constant volatility function; see \cite{HR98, BBF02} for more detailed discussions.
\par
We state two lemmas on the fair price $V(x,z)$.
\begin{lemma}
Under Assumption \ref{ubd}, it has
 $$ V_2(x)\leqslant V(x,z)\leqslant V_1(x), \quad \forall\; (x,z)\in\R_{+}^2, $$
where $V_i(x)$ is the fair price of the American put option on the time-horizon [0,T] associated with the standard GBM ,
$$ \dd X_i(t)=r X_i(t)\dd t+\sigma_i X_i(t) \dd B(t), \quad X_i(t)=x, \quad i=1,2. $$
\end{lemma}
\begin{lemma}
We have
 $$ (K-x)^+\leqslant V(x,z)\leqslant K, \quad \forall (x,z); $$
moreover, for each $z>0$, the map
\[x  \mapsto V(x,z)\]
is convex, continuous, and decreasing in $x\in [0, \infty)$.
\end{lemma}
Now we present two results. The first one is the existence of the parametric boundary $z\rightarrow b(z)$ of the option's continuation region  in $(z,x)$.  The second one
is the skewness of the time-parameter curve $t\rightarrow b(z=z(t))$ under the monotone assumption of the volatility function $\sigma(\cdot)$.

\begin{prop}\label{mainThm}
 Under Assumption \ref{ubd}, there exists a continuous curve $x=b(z)$ such that the continuation region
 $$
 C=\{(x,z)\in\R_{+}^2 \colon V(x,z)> (K-x)^+\}
 $$
 has the parametric boundary
 $$\partial C=\{ (x,z)\in\R_{+}^2 \colon x=b(z) \},$$
 and the optimal stopping time over the time-horizon $[0,T]$ is
 \begin{equation}\label{opst}
    \tau^*(\omega)=\inf\{t\in [0,T]: \ \  X(t,\omega)\leqslant b(Z(t,\omega))\}=\inf\{t\in [0,T]: \ \  X(t,\omega)= b(Z(t,\omega))\}.
 \end{equation}
 \end{prop}

 

\begin{prop}
If we assume that the volatility function $z\rightarrow\sigma (z)$ is monotone increasing in $z$, besides the   Assumption 1,
then in Proposition 3.3,  the  striking curve parametrized in the time, \/ $t\rightarrow b(z=z(t)), \ \ t \in [0,T],$\/ is not monotone increasing, though it is always squeezed by two increasing convex
curves with the same end point $(T, K).$
\end{prop}

{\bf Remark:} In \S 4.2 of \cite{HR98}, the volatility function is supposed to be 
$$\sigma(z)= \eta \sqrt{1+ \epsilon z^2}\land N,$$
 for the simulation of the volatility smile under their model. Such a  volatility function satisfies the condition of Proposition 3.4.




\section{Proofs}


{\bf Proof of Lemma 3.1.}  We use the time-change technique; see, for example, \S 5.1 of \cite{PSh06}. Define
$$
T(t,\omega )=\Big(\frac{1}{\sigma^2_2} \int_{0}^{t} \sigma^2( Z(u, \omega)) \dd u \Big) \land T, \quad t\in [0,T]. 
$$
which is strictly increasing in $t\in [0,T]$, and $T(t)\uparrow T, \ \  a.s.$ as $t\uparrow T$, by our uniform lower bound assumption on $\sigma$, namely  Assumption \ref{ubd}. The inverse
$$
\hat{T}(\theta, \omega)= \inf \{0\le t \le T: \ \ T(t,\omega )=\theta \}, \quad \theta : \theta \in [0,\big( \frac{\sigma^2_1}{\sigma^2_2}\big)T],
$$
is well-defined, and
$$
\int_{0}^{\hat{T}(\theta, \omega)} \sigma^2( Z(u, \omega)) \dd u =\theta,  \quad \theta \ge 0;
$$
moreover, $\hat{T}(\theta, \omega)$ is also  strictly increasing in $\theta$, and $\hat{T}(\theta) \uparrow T, \ \  a.s.$ as $\theta\uparrow \big( \frac{\sigma^2_1}{\sigma^2_2}\big)T$. 
Define the time-changed
motion
$$
\hat{B}(\theta, \omega)= \int_{0}^{\hat{T}(\theta, \omega)} \sigma( Z(u, \omega)) \dd B(u,\omega), \quad \theta \ge 0.
$$
Then, as  \S 5.2 \cite{PSh06} shows, the process $\theta\mapsto \hat{B}(\theta, \omega)$ is a standard Brownian motion w.r.t.
the filtration $\mathcal{\hat{F}}(\theta)= \mathcal{F}(\hat{T}(\theta)), \ \ \theta \ge 0.$
\par
Writing in them of $\theta$, we have
\begin{equation}\label{tc1}
   X(\theta)= X(0)e^{ \hat{T}(\theta)\cdot r - \frac{1}{2}\theta + \hat{B}(\theta)}.
\end{equation}
While for $X_i(\theta), \ \ i=1,2,$ we have
\begin{equation}\label{tci}
   X_i(\theta)= X(0)e^{ \frac{\theta}{\sigma_{i}^{2}}\cdot r - \frac{1}{2}\theta + B_i(\theta)},
\end{equation}
in which  $B_i(\theta))$ is the standard Brownian motion obtained from the scaling $\theta\rightarrow  \sigma_i B(\frac{\theta}{\sigma_i^2}) $. We notice that
\begin{equation}\label{thetaieq}
   \frac{\theta}{\sigma_1^2} \le   \hat{T}(\theta, \omega)        \le            \frac{\theta}{\sigma_2^2}.
\end{equation}
Since $t\rightsquigarrow \theta$ is a one-to-one transformation, we have
 $$
 V(x,z)=\esssup_{\tau'}\mathbb{E}_{(x,z)}\big[e^{-r\tau'} (K-X(\tau'))^{+}\big],
 $$
where $\tau'$ is ranging over the class of all stopping times over the scaled time-horizon  w.r.t. the time-changed Brownian filtration $\mathcal{\hat{F}}(\theta)$; so are for the
$V_i(x), \ \ i=1,2.$
\par
In term of $\theta$, $V(x,z)$ and $V_i(x)$ are all driven by the standard Brownian motion. We may compare the first term of the three exponentials in \eqref{tc1} and \eqref{tci},
together with \eqref{thetaieq}, and conclude that, for all $\theta'>0$,
$$
\mathbb{E}_{x}\big[e^{-r\theta'} (K-X_{2}(\theta'))^{+}\big]\le \mathbb{E}_{(x,z)}\big[e^{-r\theta'} (K-X(\theta'))^{+}\big] \le \mathbb{E}_{x}\big[e^{-r\theta'} (K-X_{1}(\theta'))^{+}\big].
$$
Substituting $\theta'$ by $\tau'(\theta)$, we have the desired bound.  $\Box$

\smallskip

{\bf Proof of Lemma 3.2.} Taking $\tau =0$ in the defining equality of $V(x,z)$, we have $V(x,z)\ge (K-x)^{+}$; that $V(x,z)\le K$ is obvious. Using the (\ref{dgbm}), we can
write $V(x,z)$ explicitly as
$$
V(x,z)= \esssup_{\tau}\mathbb{E}\big[e^{-r\tau'} \big(K-x\exp\{r\tau -\frac{1}{2}\int_{0}^{\tau}\sigma^2(u)\dd u+\int_{0}^{\tau}\sigma(u)\dd B(u)\}\big)^{+}\big],
$$
in which $\sigma(u)= \sigma(Z(u)), X(0)=x, Z(0)=z.$ From this display, it is seen that, for each $z$, $x\mapsto V(x,z)$ is convex, continuous, and decreasing in $x$. $\Box$

\smallskip

{\bf Proof of Proposition 3.3.} Define, for each $z$,
$$
b(z)= \sup \{ x\le K: \ \ V(x,z)= (K-x)^+\}.
$$
By Lemma 3.1,\/ $ 0< V_2(x)\leqslant V(x,z)\leqslant V_1(x) \leq K,  \forall\; (x,z)$, and the striking line for $V_i(x)$ is known to be
$x=b_i$; see, for example, \S 25.1 of \cite{PSh06},. Therefore, for each $z$,
$$
0<b_1 \leq b(z)\leq b_2 <K.
$$
Since $V(x,z)$ is continuous in $(x,z)$ (recall that $(X(t), Z(t))$ is a strong Markovian process with continuous paths, and hence it has the  Fell property), the curve $z\mapsto b(z)$ is continuous.
We claim that
\begin{equation}\label{stop}
V(x,z)= (K-x)^+, \quad \forall x: \ \ x\leq b(z);
\end{equation}
so that the curve $z\mapsto b(z)$ is indeed defining the boundary of the continuation region $C$.

Suppose, on the contrary, that, for some $0<b_0(z)<b(z),$
$$
V(b_0(z), z)> K- b_0(z).
$$
Since $V(b(z), z)= K- b(z),$ we must have, for some $\beta>1,$
$$
\frac{V(b(z), z)-V(b_0(z), z)}{b(z)-b_0(z)}= -\beta <-1.
$$
We recall that, by Lemma 3.2, $x\mapsto V(x,z)$ is decreasing.  By Lemma 3.2 again, $x\mapsto V(x,z)$ is convex, and thus we have,
$$
\frac{V(b(z), z)-V(x, z)}{b(z)-x}\leq -\beta, \quad \forall x\leq b_0(z).
$$
This will imply that
$$
V(x,z) \geq V(b(z), z) +\beta (b(z)-x) = (K- b(z)) +\beta (b(z)-x)= K + (\beta-1)b(z)- \beta x,
$$
which implies that $V(x,z)>K,$ whenever $x< (\beta-1)b(z)/ \beta.$ This is a contraction to the fact that $V(x,z)\leq K$ (Lemma 3.2). Therefore, the supposition must be false.
That the curve $z\mapsto b(z)$ lies between two parallels $x=b_i, i=1,2,$ is a consequence of Lemma 3.1.

Now we prove that the $\tau^*$ defined by (\ref{opst}) is indeed the optimal stopping time; that is,  $\tau^*(\omega)=\tau_D(\omega)$, and
$$
\tau_D(\omega)= \inf \{t\in [0, T]:\ \ (X(t,\omega), Z(t,\omega)) \in D\},
$$
where $D$ is the stopping region $D= \{(x,z): \ \ V(x,z)= (K-x)^+\}$; see \S 2.2 of \cite{PSh06}. For each $t>0$, we observe that, by the definition of $b(z)$,
$$
(X(t,\omega), Z(t,\omega)) \in D \quad \textrm{if and only if } \quad   X(t,\omega)\leq b( Z(t,\omega)).
$$
Therefore,,
\begin{eqnarray*}
  \tau_D(\omega) &=& \inf\{t\in [0, T]:\ \  X(t,\omega)\leqslant b(Z(t,\omega))\} \\
   &=& \inf\{t\in [0, T]:\ \ X(t,\omega)= b(Z(t,\omega))\} \\
   &=& \tau^*(\omega)
 \end{eqnarray*}
The second ``$=$" in the above is due to the path-continuity of the process. Indeed, suppose, on the contrary that, the ``$<$" held there, then there will exist $t_2<t_1$ such that
$$
 X(t_2,\omega) < b(Z(t_2,\omega)); \quad  X(t_1,\omega)= b(Z(t_1,\omega)).
$$
This is impossible whenever we start the process $(X,Z)$ at $(x,z)\in C$ which is above the curve $z\mapsto b(z)$.   $\Box$

\medskip

{\bf Proof of Proposition 3.4.} we claim that,
for any two $z$ and $z'$,
 $$
 (b(z)-b(z'))(V(x,z)-V(x,z'))\leqslant 0, \quad \textrm{for any $x$ between $b(z), b(z')$}.
 $$
Indeed, suppose that $b(z) <b(z')$. Then, for any $x: b(z)<x< b(z')$, by the definition of $z\mapsto b(z)$, and the proof of Theorem \ref{mainThm},
$V(x,z')= (K-x)^+$, while $V(x,z')> (K-x)^+$. Thus, $V(x,z')> V(x,z)$. On that other hand, if  $b(z) >b(z')$, then the same argument gives  $V(x,z')< V(x,z)$.

Now consider the time-parameter striking curve
 $t\rightarrow b(z(t)),\ \ t\in [0,T],$ and suppose that it is monotone increasing in $t$. Then, by the above ``anti-comonotone"  property, since $b(z)$ in increasing in $z$,  the value 
$z\rightarrow V(x,z)$ 
must be monotone decreasing, for each $x\in (\inf b(\cdot), \sup b(\cdot)).$ This is  a contraction.  We have assumed that the volatility $\sigma(z)$ is increasing in $z$,
thus as a consequence the value $V(\cdot)$ must be monotone increasing too; any option must get higher value when the volatility of the underlying asset gets higher.  Therefore the curve $t\rightarrow b(z(t))$ cannot be monotone increasing in $t\in [0,T].$ The two increasing convex curves which squeeze our striking curve are those striking curves for the two American options of each the underlying asset follows  the standard GBM
with constant volatilities $\sigma_1$ and $ \sigma_2$  respectively .  $\Box$

\smallskip

\section{Conclusion} 

1. It is the contribution of SDDE in financial economics to formulate a risky asset for which its present value brings the memory of its historical values. 
The choice of the memory $Y(t)=X(t-t_0)$ mentioned in Section 1 reflects that a past time-instant $t_0$ is the source of the ``after-effect", and this can be extended to $n$ time-instants
$t_0,\cdots, t_n$. The European option pricing  based on the  SDDE of this type appeared in \cite{AHMP07}, \cite{KSW07}, and \cite{MS13}; these papers  rely on the martingale aspect of the pricing theory.
The American option pricing is  certainly at the beginning to be viewed from the martingale aspect, as one may see from \cite{PSh06}. However it is more important then to move to view American options  in the Markov process aspect; since only then the striking curve can be discussed, that is, the parametric curve to separate the region in which the owner of the option holds and waits, and the region in which the owner exercises and gets the (positive) reward.  In the basic ( that is the constant volatility is assumed) BSM theory, the striking curve
is a monotone increasing and  convex curve across the time horizon $[0,T]$; see  Chapter 8 of \cite{S04} or Section 25.1 of \cite{PSh06}.
The novelty of this article  is  that, if we assume the volatility is the ratio of the asset's
present value $X(t)$ and historical value $Y(t)$, with the choice of the historical values being  exponentially averaged, then a parametric striking curve still appears, \emph{yet} it is
a ``anti-comonotone" curve,  as shown in Proposition 3.4. This would assert  that the striking curve is \emph{skewed}, due to the historical value of the asset. We mention that, to our knowledge, this situation is  firstly observed, and we would compare 
this result  with  one main conclusion  in \cite{HR98},  in which the authors  discuss the volatility smile of European options under the model.

2. Financial economics under  uncertainly is one  fundamental topic in Microeconomics, and we refer to Chapter 6 of \cite{MWG95}. Option pricing is one aspect, and here we would discuss
of the effect of the asset's historical value to the pricing turnout; the volatility smile for European options in \cite{HR98}, and the striking-curve skewness  for American options in this article.    Portfolio selection (under uncertainty) is also one classical aspect, it can be traced to the classics \cite{NM44}, and we  cite two very recent papers in this aspect \cite{CPZR17} and \cite{QKZ16}. Study of  portfolio selections of risky assets with the memory effects seems to be very promising.




\vspace{.2in}

\begin{quote}
\begin{small}

\noindent \textsc{Narn-Rueih Shieh}.\ Mathematics Department,
National
Taiwan University, Taipei 10617, Taiwan.
 E-mail: \texttt{shiehnr@ntu.edu.tw } \\


\end{small}
\end{quote}


\end{document}